\documentclass{article}\usepackage{amsmath,amssymb,latexsym,amsfonts,enumerate,theorem,longtable,amscd}
\newtheorem{definition}{Definition}[section]\newtheorem{theorem}[definition]{Theorem}\newtheorem{proposition}[definition]{Proposition}
\newtheorem{lemma}[definition]{Lemma}\newtheorem{example}{Example}
\begin{document}

\begin{center}
{\large\bf On $q$- and $h$-deformations of 3d-superspace}
\end{center}

\begin{center}
Salih Celik

Department of Mathematics, Yildiz Technical University, DAVUTPASA-Esenler, Istanbul, 34220 TURKEY.
\end{center}



\begin{abstract}
In this paper, we introduce non-standard deformations of (1+2)- and (2+1)-superspaces via a contraction using standard deformations of them. This deformed superspaces denoted by ${\mathbb A}_h^{1|2}$ and ${\mathbb A}_{h'}^{2|1}$, respectively. We find a two-parameter $R$-matrix satisfying quantum Yang-Baxter equation and thus obtain a {\it new} two-parameter non-standard deformation of the supergroup ${\rm GL}(1|2)$. Finally, we get a new superalgebra derived from the super-Hopf algebra of functions on the quantum superspace ${\mathbb A}_{p,q}^{1|2}$.
\end{abstract}

\section{Introduction}\label{sec1}
There are two distinct deformations for general Lie (super)groups as standard and non-standard (or Jordanian). One of them is the well-known quantum ($q$-deformed) group and the other is the so-called Jordanian ($h$-deformed) one. Specially, quantum groups ${\rm GL}_q(2)$ \cite{Manin1} and ${\rm GL}_h(2)$ \cite{kup} have been obtained by deforming the coordinates of a plane to be noncommutative objects. In \cite{agha}, the authors have shown that the $h$-deformed group can be obtained from the $q$-deformed Lie group through a singular limit $q\to1$ of a linear transformation. This method is known as the contraction procedure. Using this method, one- and two-parameter $h$-deformations of supergroup ${\rm GL}(1|1)$ were obtained in \cite{dab} and \cite{celik1}, respectively.

In this paper, we give some standard (as $q$-deformation) deformations of (1+2)-superspace using the super-Hopf algebra structure of ${\cal O}({\mathbb A}^{1|2})$ and non-standard (as $h$-deformation) deformations using standard deformations via a contraction. We also introduce an $(h,h')$-deformed supergroup acting these two parameter $h$-deformed superspaces. Finally, we define involutions on $h$-deformed superspaces and use the generators of $(p,q)$-deformed super algebra ${\cal O}({\mathbb A}_{p,q}^{1|2})$ to get a new Lie superalgebra.

Throughout the paper, we will fix a base field ${\mathbb K}$. The reader may consider it as the set of real numbers, ${\mathbb R}$, or the set of complex numbers, ${\mathbb C}$. We will denote by ${\mathbb G}$ the Grassmann numbers and by ${\mathbb K}'$ the set ${\mathbb K}\cup{\mathbb G}$.

\section{On $(p,q)$-deformation of superspaces ${\mathbb A}^{1|2}$ and ${\mathbb A}^{2|1}$}\label{sec2}

In order to define superalgebras and Hopf superalgebras, one makes some minor changes in familiar definitions. For example, if a space ${\cal X}$ is a superspace or ${\mathbb Z}_2$-graded, then we denote by $\tau(a)$ the ${\mathbb Z}_2$-degree of the element $a\in{\cal X}$. If $\tau(a)=0$, then we will call the element $a$ {\it even} and if $\tau(a)=1$, it is called {\it odd}.


\subsection{The algebra of polynomials on the quantum superspace ${\mathbb A}^{1|2}_q$}\label{sec2.1}
Let ${\mathbb K}\langle X,\Theta_1,\Theta_2 \rangle$ be a free algebra with unit generated by $X$, $\Theta_1$ and $\Theta_2$, where the coordinate $X$ is even, the coordinates $\Theta_1$ and $\Theta_2$ are odd.

\begin{definition} \cite{Manin2} 
Let $I_q$ be the two-sided ideal of ${\mathbb K}\langle X,\Theta_1,\Theta_2 \rangle$ which generated by the elements $X\Theta_1-q\Theta_1 X$, $X\Theta_2-q\Theta_2 X$, $\Theta_1\Theta_2+q^{-1}\Theta_2\Theta_1$, $\Theta_1^2$ and $\Theta_2^2$. The quantum superspace ${\mathbb A}^{1|2}_q$ with the function algebra
$${\cal O}({\mathbb A}^{1|2}_q) = {\mathbb K}\langle X,\Theta_1,\Theta_2 \rangle/I_q$$
is called ${\mathbb Z}_2$-graded quantum space (or quantum superspace).
\end{definition}
This associative algebra over the complex number is known as the algebra of polynomials over quantum (1+2)-superspace. In accordance with Definition 2.1, we have
\begin{equation} \label{2.1}
X\Theta_i = q\Theta_i X, \quad \Theta_i\Theta_j = -q^{i-j}\Theta_j\Theta_i, \qquad (i,j=1,2)
\end{equation}
where $q\in{\mathbb K}-\{0\}$.

\begin{example} 
If we consider the generators of the algebra ${\cal O}({\mathbb A}^{1|2}_q)$ as linear maps, then we can find the matrix representations of them. In fact,
it can be seen that there exists a representation $\rho:{\cal O}({\mathbb A}^{1|2}_q)\to \mbox{M}(3,{\mathbb K}')$ such that matrices
\begin{equation} \label{2.2}
\rho(X)=\begin{pmatrix} q & 0 & 0 \\ 0 & q & 0 \\ 0 & 0 & q^2 \end{pmatrix}, \quad
\rho(\Theta_1)=\begin{pmatrix} 0 & 0 & \varepsilon_1 \\ 0 & 0 & 0 \\ 0 & 0 & 0 \end{pmatrix}, \quad
\rho(\Theta_2)=\begin{pmatrix} 0 & 0 & 0 \\ 0 & 0 & \varepsilon_2 \\ 0 & 0 & 0 \end{pmatrix}
\end{equation}
representing the coordinate functions satisfy relations $(\ref{2.1})$ for all $\varepsilon_1, \varepsilon_2$.
\end{example}


\noindent{\bf Note 1.} In the next section, we will assume that $\varepsilon_1$ and $\varepsilon_2$ are two Grassmann numbers.

The following definition gives the product rule for tensor product of ${\mathbb Z}_2$-graded algebras.

\begin{definition} 
The product rule is defined by
\begin{equation*}
(a_1\otimes a_2)(a_3\otimes a_4)=(-1)^{\tau(a_2)\tau(a_3)}(a_1a_3\otimes a_2a_4)
\end{equation*}
in the ${\mathbb Z}_2$-graded algebra ${\cal A}\otimes{\cal A}$, where ${\cal A}$ is the ${\mathbb Z}_2$-graded algebra and $a_i$'s are homogeneous elements in ${\cal A}$.
\end{definition}

We denote the unital extension of ${\cal O}({\mathbb A}^{1|2}_q)$ by ${\cal F}({\mathbb A}^{1|2}_q)$ adding the unit and $x^{-1}$, the inverse of $x$, which obeys $x x^{-1} = {\bf 1} =  x^{-1} x$. The following theorem says that the superalgebra ${\cal F}({\mathbb A}_q^{1|2})$ has a Hopf algebra structure \cite{celik2}.

\begin{theorem} 
The superalgebra ${\cal F}({\mathbb A}_q^{1|2})$ is a super-Hopf algebra with the defining coproduct, counit and coinverse on the algebra
${\cal F}({\mathbb A}_q^{1|2})$ as follows:

\noindent$(1)$ The coproduct $\Delta: {\cal F}({\mathbb A}_q^{1|2}) \longrightarrow {\cal F}({\mathbb A}_q^{1|2}) \otimes {\cal F}({\mathbb A}_q^{1|2})$ is defined by
\begin{equation} \label{2.3}
\Delta(X) = X \otimes X, \quad \Delta(\Theta_i) = \Theta_i \otimes X^i + X^i \otimes \Theta_i, \qquad (i=1,2).
\end{equation}

\noindent$(2)$ The counit $\epsilon: {\cal F}({\mathbb A}_q^{1|2}) \longrightarrow {\mathbb K}$ is given by
\begin{equation*}
\epsilon(X)=1, \quad \epsilon(\Theta_i)=0, \qquad (i=1,2).
\end{equation*}

\noindent$(3)$ The algebra ${\cal F}({\mathbb A}_q^{1|2})$ admits a ${\mathbb K}$-algebra antihomomorphism (coinverse)
$S:{\cal F}({\mathbb A}_q^{1|2}) \longrightarrow {\cal F}({\mathbb A}^{1|2}_{q^{-1}})$ defined by
\begin{equation*}
S(X) = X^{-1}, \quad S(\Theta_i) = - X^{-i} \Theta_i X^{-i}, \qquad (i=1,2).
\end{equation*}
\end{theorem}

\subsection{The algebra of polynomials on the quantum superspace ${\mathbb A}^{2|1}_{p,q}$}\label{sec2.2}

Let ${\mathbb K}\langle \Phi,Y_1,Y_2 \rangle$ be a free algebra with unit generated by $\Phi$, $Y_1$ and $Y_2$, where $\tau(\Phi)=1$, and $\tau(Y_1)=0=\tau(Y_2)$.

\begin{definition} \cite{sucelik1} 
Let $\Lambda({\mathbb A}_q^{1|2})$ be the algebra with the generators $\Phi$, $Y_1$ and $Y_2$ satisfying the relations
\begin{equation} \label{2.4}
\Phi^2 = 0, \quad \Phi Y_1 = qp^{-1} Y_1 \Phi, \quad \Phi Y_2 = pq Y_2 \Phi, \quad Y_1Y_2 = pq^{-1} \,Y_2Y_1.
\end{equation}
We call $\Lambda({\mathbb A}_q^{1|2})$ exterior algebra of the ${\mathbb Z}_2$-graded space ${\mathbb A}_q^{1|2}$.
\end{definition}

\noindent{\bf Note 2.} The exterior algebra $\Lambda({\mathbb A}_q^{1|2})$ of the superspace ${\mathbb A}_q^{1|2}$ can be thought of as a two-parameter deformation of the (2+1)-superspace ${\mathbb A}^{2|1}$. So, we denote this algebra by ${\cal O}({\mathbb A}^{2|1}_{p,q})$.

\begin{example} 
If we consider the generators of the algebra ${\cal O}({\mathbb A}^{2|1}_{p,q})$ as linear maps, then we can find the matrix representations of them. In fact,
it can be seen that there exists a representation $\rho:{\cal O}({\mathbb A}^{2|1}_{p,q})\to \mbox{M}(3,{\mathbb K}')$ such that matrices
\begin{equation*}
\rho(\Phi)=\begin{pmatrix} 0 & 0 & \epsilon \\ 0 & 0 & 0 \\ 0 & 0 & 0 \end{pmatrix}, \quad
\rho(Y_1)=\begin{pmatrix} q & 0 & 0 \\ 0 & p & 0 \\ 0 & 0 & p \end{pmatrix}, \quad
\rho(Y_2)=\begin{pmatrix} 0 & c & 0 \\ 0 & 0 & 0 \\ 0 & 0 & 0 \end{pmatrix}
\end{equation*}
representing the coordinate functions satisfy relations $(\ref{2.4})$ for all $c, \varepsilon$.
\end{example}

\section{Two parameter $h$-deformation of the superspaces} \label{sec3}

In this section, we introduce a two-parameter $h$-deformation of the superspace ${\mathbb A}^{1|2}$ (and its dual) from the $(p,q)$-deformation via a contraction similar to the method of \cite{agha}.

We consider the $q$-deformed algebra of functions on the quantum superspace ${\mathbb A}_q^{1|2}$ generated by $X$, $\Theta_1$ and $\Theta_2$ with the relations (\ref{2.1}) and we introduce new even coordinate $x$ and odd coordinates $\theta_1$, $\theta_2$ with the change of basis in the coordinates of the $q$-superspace using the following $g$ matrix:
\begin{equation} \label{3.1}
{\bf X} = \begin{pmatrix} X \\ \Theta_1 \\ \Theta_2 \end{pmatrix} = \begin{pmatrix} 1 & 0 & \tilde{h}' \\ 0 & 1 & 0 \\ \tilde{h} & 0 & 1 \end{pmatrix} \begin{pmatrix} x \\ \theta_1 \\ \theta_2 \end{pmatrix} = g \,{\bf x}, \quad \tilde{h}=\frac{h}{q-1}, \quad \tilde{h}'=\frac{h'}{pq-1}
\end{equation}
where $h$ and $h'$ $(h\ne0\ne h')$ are two {\it new} deformation parameters that will be replaced with $q$ and $p$ $(q\ne1\ne pq)$ in the limits $q\to1$ and $p\to1$.

We now assume that the parameters $h$ and $h'$ are both Grassmann numbers ($h^2=0=h'^2, \,\, hh'=-h'h$) and anticommute with $\theta_i$ for $i=1,2$. When the relations (\ref{2.1}) are used, one gets
\begin{equation} \label{3.2}
x\theta_1 = q\theta_1 x, \,\, x\theta_2 = q\theta_2 x + h x^2, \,\, \theta_2\theta_1 = -q\theta_1\theta_2, \,\, \theta_1^2 = 0, \,\,
\theta_2^2 = -h\theta_2 x.
\end{equation}
Note that the parameter $h'$ does not enter the above relations. By taking the limit $q\to1$ we obtain the following exchange relations, which define the $h$-superspace ${\mathbb A}_{h}^{1|2}$:

\begin{definition} \cite{celik2} 
Let ${\cal O}({\mathbb A}_{h}^{1|2})$ be the algebra with the generators $x$, $\theta_1$ and $\theta_2$ satisfying the relations
\begin{equation} \label{3.3}
x\theta_1 = \theta_1 x, \quad x\theta_2 = \theta_2 x + h x^2, \quad \theta_1\theta_2 = -\theta_2\theta_1, \quad \theta_1^2 = 0, \quad
\theta_2^2 = -h\theta_2 x.
\end{equation}
We call ${\cal O}({\mathbb A}_{h}^{1|2})$ the algebra of functions on the ${\mathbb Z}_2$-graded quantum space ${\mathbb A}_{h}^{1|2}$.
\end{definition}

\begin{example} 
Let us assume that $\varepsilon_1$ and $\varepsilon_2$ are two Grassmann numbers. If the $g$ matrix in $(\ref{3.1})$ is used, the matrix representation in $(\ref{2.2})$ takes the following form:
\begin{equation*}
\rho(x)=\begin{pmatrix} q(1-\tilde{h}\tilde{h}') & 0 & 0 \\ 0 & q(1-\tilde{h}\tilde{h}') & -\tilde{h}'\varepsilon_2 \\ 0 & 0 & q^2(1-\tilde{h}\tilde{h}') \end{pmatrix}, \quad \rho(\theta_1)=\begin{pmatrix} 0 & 0 & \varepsilon_1 \\ 0 & 0 & 0 \\ 0 & 0 & 0 \end{pmatrix},
\end{equation*}
\begin{equation} \label{3.4}
\rho(\theta_2)=\begin{pmatrix} -q\tilde{h} & 0 & 0 \\ 0 & -q\tilde{h} & (1+\tilde{h}\tilde{h}')\varepsilon_2 \\ 0 & 0 & -q^2\tilde{h} \end{pmatrix}.
\end{equation}
These matrices satisfy the relations $(\ref{3.2})$, for all $\varepsilon_1$ and $\varepsilon_2$.
\end{example}

\noindent{\it Proof} 
Existing claims come from the fact that $\rho$ is an algebra homomorphism. $\square$

In the case of dual (exterior) $h'$-superspace, 
we use the transformation
\begin{eqnarray} \label{3.5}
\hat{\bf X}=g\hat{\bf x}
\end{eqnarray}
with the components $\varphi$, $y_1$ and $y_2$ of $\hat{\bf x}$. The definition is given below.

\begin{definition} 
Let ${\cal O}({\mathbb A}_{h'}^{2|1}):=\Lambda({\mathbb A}_{h}^{1|2})$ be the algebra with the generators $\varphi$, $y_1$ and $y_2$ satisfying the relations
\begin{eqnarray} \label{3.6}
\varphi y_1 = y_1 \varphi, \quad \varphi y_2 = y_2 \varphi + h' y_2^2, \quad y_1y_2 = y_2y_1, \quad \varphi^2 = h'y_2\varphi
\end{eqnarray}
where $\tau(\varphi)=1$ and $\tau(y_1)=0=\tau(y_2)$. We call $\Lambda({\mathbb A}_{h}^{1|2})$ the quantum exterior algebra of the ${\mathbb Z}_2$-graded quantum space ${\mathbb A}_{h}^{1|2}$.
\end{definition}

\noindent{\bf Note 3.} The parameter $h$ does not enter the relations (\ref{3.6}). The exterior algebra $\Lambda({\mathbb A}_{h}^{1|2})$ of the superspace ${\mathbb A}_{h}^{1|2}$ can be thought of as an $h'$-deformation of the (2+1)-superspace ${\mathbb A}^{2|1}$.

\section{An $R$-matrix and its some properties} \label{sec4}
The relations in (\ref{2.1}) can be written in a compact form as follows:
\begin{equation} \label{4.1}
p \, {\bf X}\otimes{\bf X} = \hat{R}_{p,q} \,{\bf X}\otimes{\bf X}
\end{equation}
with an $R$-matrix given by \cite{sucelik2}
\begin{equation*}
\hat{R}_{p,q} = \begin{pmatrix}
p & 0 & 0 & 0 & 0 & 0 & 0 & 0 & 0 \\
0 & p-1 & 0 & q & 0 & 0 & 0 & 0 & 0 \\
0 & 0 & 0 & 0 & 0 & 0 & pq & 0 & 0 \\
0 & pq^{-1} & 0 & 0 & 0 & 0 & 0 & 0 & 0 \\
0 & 0 & 0 & 0 & -1 & 0 & 0 & 0 & 0 \\
0 & 0 & 0 & 0 & 0 & 0 & 0 & -pq^{-1} & 0 \\
0 & 0 & q^{-1} & 0 & 0 & 0 & p-1 & 0 & 0\\
0 & 0 & 0 & 0 & 0 & -q & 0 & p-1 & 0 \\
0 & 0 & 0 & 0 & 0 & 0 & 0 & 0 & -1 \\
\end{pmatrix}
\end{equation*}
where $p,q\in{\mathbb K}-\{0\}$. This matrix satisfies the graded braid equation and the matrix $R_{p,q}=P\hat{R}_{p,q}$ satisfies the graded Yang-Baxter equation where $P$ is the super permutation matrix.

It can be considered that a change of basis in the quantum superspaces leads to a two-parameter $R$-matrix. The corresponding $R$-matrix can be obtained as
\begin{equation*} \label{4.2}
\hat{R}_{h,h'} = \lim_{(p,q)\to(1,1)} \left[(g\otimes g)^{-1} \hat{R}_{p,q}(g\otimes g)\right]
\end{equation*}
where it is assumed that $\otimes$ is graded. As a result, we obtain the following $R$-matrix
\begin{equation*} \label{4.3}
\hat{R}_{h,h'} = \begin{pmatrix}
1+hh' & 0 & h' & 0 & 0 & 0 & -h' & 0 & 0 \\
0 & 0 & 0 & 1 & 0 & 0 & 0 & 0 & 0 \\
h & 0 & hh' & 0 & 0 & 0 & 1 & 0 & -h' \\
0 & 1 & 0 & 0 & 0 & 0 & 0 & 0 & 0 \\
0 & 0 & 0 & 0 & -1 & 0 & 0 & 0 & 0 \\
0 & 0 & 0 & 0 & 0 & 0 & 0 & -1 & 0 \\
-h & 0 & 1 & 0 & 0 & 0 & hh' & 0 & -h'\\
0 & 0 & 0 & 0 & 0 & -1 & 0 & 0 & 0 \\
0 & 0 & -h & 0 & 0 & 0 & -h & 0 & hh'-1 \\
\end{pmatrix}.
\end{equation*}

The equation in (\ref{4.1}) with the {\it new} $R$-matrix $\hat{R}_{h,h'}$ takes the form
\begin{equation*} \label{4.4}
{\bf x}\otimes{\bf x} = \hat{R}_{h,h'} \,{\bf x}\otimes{\bf x},
\end{equation*}
that is, the relations (\ref{3.3}) are equivalent to this equation.

The $R$-matrix $\hat{R}_{h,h'}$ has some interesting properties. Some of them are listed below, where sometimes we write $\hat{R}=\hat{R}_{h,h'}$ for simplicity.
\begin{enumerate}
\item The matrix $\hat{R}_{h,h'}$ satisfies the graded (and ungraded) braid equation $\hat{R}_{12}\hat{R}_{23}\hat{R}_{12} = \hat{R}_{23}\hat{R}_{12}\hat{R}_{23}$, where $\hat{R}_{12}=\hat{R} \otimes I_3$ and $\hat{R}_{12}=I_3 \otimes \hat{R}$.
\item The matrix $R_{h,h'}=P\hat{R}_{h,h'}$ satisfies the graded (and ungraded) Yang-Baxter equation $R_{12}R_{13}R_{23} = R_{23}R_{13}R_{12}$,
    where $R_{13}$ acts both on the first and third spaces (both in graded and ungraded).
\item The matrix $\hat{R}_{h,h'}$ holds $\hat{R}_{h,h'}^2=I_9$ and thus it has two eigenvalues $\pm1$.
\item If we set $hh'=0$, then the matrix $R_{h,h'}$ can be decomposed in the form
\begin{equation*}
R_{h,h'} = R(h) R(h')
\end{equation*}
where
\begin{equation*}
R(h) = \begin{pmatrix}
1 & 0 & 0 & 0 & 0 & 0 & 0 & 0 & 0 \\
0 & 1 & 0 & 0 & 0 & 0 & 0 & 0 & 0 \\
-h & 0 & 1 & 0 & 0 & 0 & 0 & 0 & 0 \\
0 & 0 & 0 & 1 & 0 & 0 & 0 & 0 & 0 \\
0 & 0 & 0 & 0 & 1 & 0 & 0 & 0 & 0 \\
0 & 0 & 0 & 0 & 0 & 1 & 0 & 0 & 0 \\
h & 0 & 1 & 0 & 0 & 0 & 1 & 0 & 0\\
0 & 0 & 0 & 0 & 0 & 0 & 0 & 1 & 0 \\
0 & 0 & h & 0 & 0 & 0 & h & 0 & 1 \\
\end{pmatrix}, \quad R(h') = R^{\rm st}(h)|_{h=h'}.
\end{equation*}
   It can be checked that these matrices both satisfy the graded (and ungraded) Yang-Baxter equation.
\item If $P_\pm$ are the projections onto the eigenspaces $\pm1$ of $\hat{R}_{h,h'}$, then we have
\begin{equation*}
\hat{R}_{h,h'} = P_+ - P_-.
\end{equation*}
Let ${\cal O}({\mathbb A}^{1|2})$ and ${\cal O}({\mathbb A}^{2|1})$ be the quotients of algebras generated by $x$, $\theta_1$, $\theta_2$ and $\varphi$, $y_1$, $y_2$ modulo the two-sided ideals generated by ${\rm Ker}P_-$ and ${\rm Ker}P_+$, respectively. Then ${\cal O}({\mathbb A}^{1|2})$ and ${\cal O}({\mathbb A}^{2|1})$ are isomorphic to ${\cal O}({\mathbb A}_{h}^{1|2})$ with defining relations (\ref{3.3}) and ${\cal O}({\mathbb A}_{h'}^{2|1})$ with defining relations (\ref{3.6}), respectively. That is, we can write
\begin{equation*}
P_- \, {\bf x}\otimes{\bf x} = 0 \quad {\rm and} \quad (-1)^{\tau(\hat{\bf x})} \,P_+ \,\hat{\bf x}\otimes\hat{\bf x} = 0.
\end{equation*}
\end{enumerate}

\section{The quantum super bialgebra ${\cal O}({\rm M}_{h,h'}(1|2))$}\label{sec5}

Let $T$ be a 3x3 matrix in ${\mathbb Z}_2$-graded space given by
\begin{equation*}
T = \begin{pmatrix} a & \alpha & \beta \\ \gamma & b & c \\ \delta & d & e \end{pmatrix} =(t_{ij})
\end{equation*}
where $a$, $b$, $c$, $d$, $e$ are even and  $\alpha$, $\beta$, $\gamma$ and $\delta$ are odd. The coordinate ring of such matrices over a field ${\mathbb K}$ is simply the polynomial ring in nine variables, that is ${\cal O}({\rm M}(1|2))={\mathbb K}[a,b,c,d,e,\alpha,\beta,\gamma,\delta]$.

In this section, we will assume that the matrix entries of $T$ belong a free superalgebra and define a two-parameter $h$-analogue of ${\cal O}({\rm M}(1|2))$. To do so, let $x$, $\theta_1$, $\theta_2$ be elements of the superalgebra ${\cal O}({\mathbb A}_{h}^{1|2})$ subject to the relations (\ref{3.3}) and $\varphi$, $y_1$, $y_2$ be elements of ${\cal O}({\mathbb A}_{h'}^{2|1})$ subject to the relations (\ref{3.6}), and $t_{ij}$ be nine generators which supercommute with the elements of ${\cal O}({\mathbb A}_{h}^{1|2})$ and ${\cal O}({\mathbb A}_{h'}^{2|1})$. It is well known that, the supermatrix $T$ defines the linear transformations $T:{\mathbb A}_{h}^{1|2}\longrightarrow {\mathbb A}_{h}^{1|2}$ and $T:{\mathbb A}_{h'}^{2|1}\longrightarrow {\mathbb A}_{h'}^{2|1}$. Let
${\bf x}=(x,\theta_1,\theta_2)^t$ and $\hat{\bf x}=(\varphi,y_1,y_2)^t$. So, we can give the following theorem.

\begin{theorem} 
Under the above hypotheses, the following conditions are \\ equivalent:

\noindent$(i)$ $T{\bf x}={\bf x}'\in{\mathbb A}_{h}^{1|2}$ and $T\hat{\bf x}=\hat{\bf x}'\in{\mathbb A}_{h'}^{2|1}$,

\noindent$(ii)$ the relations are satisfied
\begin{align} \label{5.1}
a\alpha &= (1+hh')\alpha a - h'(\alpha\delta + da), \quad a\beta = \beta a + h'(a^2 - ea - \beta\delta) - h\beta^2, \nonumber\\
a\gamma &= (1+hh')\gamma a + h(\gamma\beta - ca), \quad ac = ca - hc\beta - h'\gamma a + hh'\gamma\beta,\nonumber\\
a\delta &= \delta a + h(a^2 -ea + \delta\beta) + h'\delta^2, \quad ad = da + h\alpha a + h'd\delta - hh'\alpha\delta, \nonumber\\
ae &= ea + h\beta(a-e) + h'(e-a)\delta, \quad \alpha\beta=-(1+hh')\beta\alpha + h'(\beta d + e\alpha), \nonumber\\
\alpha\gamma &= -\gamma\alpha, \quad \alpha c = c\alpha, \quad \alpha\delta = -\delta\alpha - ha\alpha + h' \delta d - hh' ad, \nonumber\\
\alpha d &= d\alpha + h'd^2, \quad \alpha e = e\alpha + h\beta\alpha + h' ed - hh' d\beta, \nonumber\\
\beta\gamma &= - \gamma\beta + hc\beta - h'\gamma a - hh' ca, \quad \beta c = (1-hh')c\beta - h'(\gamma\beta + ca), \nonumber\\
\beta\delta &= - \delta\beta + (h\beta + h'\delta)(e - a), \quad \beta d = d\beta + h\alpha\beta + h'de - hh' e\alpha, \nonumber\\
\beta e &= e\beta + h'(e^2 - ea - \delta\beta) - h\beta^2, \quad \gamma c = c\gamma + hc^2, \nonumber\\
\gamma\delta &= - (1+hh')\delta\gamma + h(e\gamma + \delta c), \quad \gamma d = d\gamma, \\
\gamma e &= e\gamma + hec - h'\delta\gamma - hh' c\delta, \quad c\delta = \delta c - h ec - h'\delta\gamma - hh' \gamma e, \nonumber\\
cd &= dc, \quad ce = (1-hh')ec + h'(e\gamma - \delta c), \,\,\, \delta d = (1-hh')d\delta + h(\alpha\delta - da), \nonumber\\
\delta e &= e\delta + h(e^2 - ea + \beta\delta) + h'\delta^2, \quad de = (1-hh')ed + h(\beta d - e\alpha), \nonumber\\
\alpha^2 &= h'\alpha d, \quad \beta^2 = h'\beta (e-a), \quad \gamma^2 = h\gamma c, \quad \delta^2 = h\delta (e-a), \nonumber\\
b \,t_{ij} &= t_{ij} b, \quad a(h\beta + h'\delta) = (h\beta + h'\delta)a, \quad e(h\beta + h'\delta) = (h\beta + h'\delta)e. \nonumber
\end{align}
\end{theorem}

\noindent{\it Proof} 
A direct verification shows that the relations (\ref{5.1}) respect the ideals defining ${\mathbb A}_{h}^{1|2}$ and ${\mathbb A}_{h'}^{2|1}$. $\square$

Standard FRT construction \cite{FRT}, namely, the relations (\ref{5.1}), is obtained via the matrix $\hat{R}_{h,h'}$ given in Sect. 4:

\begin{theorem} 
A 3x3-matrix $T$ is a ${\mathbb Z}_2$-graded quantum supermatrix if and only if
\begin{equation*}
\hat{R}_{h,h'} T_1T_2 = T_1T_2 \hat{R}_{h,h'}
\end{equation*}
where $T_1=T\otimes I_3$ and $T_2=PT_1P$.
\end{theorem}

\begin{definition} 
The superalgebra ${\cal O}({\rm M}_{h,h'}(1|2))$ is the quotient of the free algebra ${\mathbb K}\{a,b,c,d,e,\alpha,\beta,\gamma,\delta\}$ by the two-sided ideal $J_{h,h'}$ generated by the relations $(\ref{5.1})$ of Theorem 5.1.
\end{definition}

\noindent{\bf Note 4.} The quantum matrix space ${\rm M}_{p,q}(1|2)$ is obtained in \cite{sucelik2}. It is clear that a change of basis in the quantum superspace leads to the similarity transformation $T = g^{-1}T'g$, where $T'\in {\rm M}_{p,q}(1|2)$. Therefore, the entries of the transformed quantum matrix $T$ fulfill the commutation relations (\ref{5.1}) of the matrix elements of the matrix $T$ in ${\rm M}(1|2)$.

\begin{theorem} 
The superalgebra ${\cal O}({\rm M}_{h,h'}(1|2))$ with the following two algebra homomorphisms of superalgebras

\noindent$(1)$ the coproduct $\Delta: {\cal O}({\rm M}_{h,h'}(1|2))\longrightarrow {\cal O}({\rm M}_{h,h'}(1|2))\otimes{\cal O}({\rm M}_{h,h'}(1|2))$ determined by $\Delta(t_{ij})=\sum_{k=1}^3 t_{ik}\otimes t_{kj}$,

\noindent$(2)$ the counit $\epsilon: {\cal O}({\rm M}_{h,h'}(1|2))\longrightarrow {\mathbb K}$ determined by $\epsilon(t_{ij})=\delta_{ij}$

\noindent becomes a super bialgebra.
\end{theorem}

\noindent{\it Proof} 
It can be easily checked the properties of the costructures hold:

\noindent(i) The coproduct $\Delta$ is coassociative in the sense of
\begin{equation*}
(\Delta\otimes{\rm id})\circ\Delta = ({\rm id}\otimes\Delta)\circ\Delta
\end{equation*}
where id denotes the identity map on ${\rm M}_{h,h'}(1|2)$ and $\Delta(ab)=\Delta(a) \Delta(b)$, $\Delta({\bf 1})={\bf 1}\otimes {\bf 1}$.

\noindent(ii) The counit $\epsilon$ has the property
\begin{equation*}
m\circ(\epsilon\otimes{\rm id})\circ\Delta = {\rm id} = m\circ({\rm id}\otimes\epsilon)\circ\Delta
\end{equation*}
where $m$ stands for the algebra product and $\epsilon(ab)=\epsilon(a) \epsilon(b)$, $\epsilon({\bf 1})=1$. $\square$

It is well known that ${\cal O}({\mathbb A}^{1|2})$ is comodule algebra over the bialgebra ${\cal O}({\rm M}(1|2))$. The following
theorem gives a quantum version of this fact.

\begin{theorem} 
There exist algebra homomorphisms
\begin{align*} \label{5.2}
\delta_L:{\cal O}({\mathbb A}_{h}^{1|2})\longrightarrow{\cal O}({\rm M}_{h,h'}(1|2))\otimes {\cal O}({\mathbb A}_{h}^{1|2}), \quad
\delta_L(x_i) = \sum_{k=1}^3 t_{ik} \otimes x_k, \nonumber\\
\tilde{\delta}_L:{\cal O}({\mathbb A}_{h'}^{2|1})\longrightarrow{\cal O}({\rm M}_{h,h'}(1|2))\otimes {\cal O}({\mathbb A}_{h'}^{2|1}), \quad
\tilde{\delta}_L(\hat{x}_i) = \sum_{k=1}^3 t_{ik} \otimes \hat{x}_k
\end{align*}
where $x_i \in \{x, \theta_1, \theta_2\}$ and $\hat{x}_i \in \{\varphi, y_1, y_2\}$.
\end{theorem}

\noindent{\it Proof} 
Using the relations (\ref{3.3}) and (\ref{3.6}) together with (\ref{5.1}), it is enough to check that
\begin{equation*}
\delta_L(x\theta_1-\theta_1 x)=\delta_L(x)\delta_L(\theta_1)-\delta_L(\theta_1)\delta_L(x)=0,
\end{equation*}
etc., in ${\cal O}({\rm M}_{h,h'}(1|2))\otimes {\cal O}({\mathbb A}_{h}^{1|2})$. To see that $\delta_L$ defines a comodule structure we check that
\begin{equation*}
(\Delta\otimes{\rm id})\circ\delta_L = ({\rm id}\otimes\delta_L)\circ\delta_L, \quad m\circ(\epsilon\otimes{\rm id})\circ\delta_L = {\rm id}. \hfill \square
\end{equation*}

A quantum supergroup (super-Hopf algebra) can be regarded as a generalization of the notion of a supergroup. It is defined by
\begin{equation*}
{\cal O}({\rm GL}_{h,h'}(1|2)) = {\cal O}({\rm M}_{h,h'}(1|2))[t]/(t \,{\rm sdet}_{h,h'}-1).
\end{equation*}
This case is also inviting to also generalize the corresponding notions of differential geometry \cite{woronowicz}. A differential calculus on
${\cal O}({\rm GL}_{h,h'}(1|2))$ will be discussed in the next work.

\section{A new superalgebra derived from ${\cal F}({\mathbb A}_{p,q}^{1|2})$}\label{sec6}

It is known that an element of a Lie group can be represented by exponential of an element of its Lie algebra. In \cite{celik3}, by virtue of this fact, using the generators of the superalgebra ${\cal F}({\mathbb A}_q^{1|1})$ it has been obtained a new super-algebra from this algebra. In this section, we will obtain a new super-algebra from ${\cal F}({\mathbb A}_{p,q}^{1|2})$. So, let us begin with definition of ${\cal F}({\mathbb A}_{p,q}^{1|2})$ which is an extension to two parameter of ${\cal F}({\mathbb A}_q^{1|2})$.

\begin{definition} 
Let $I_{p,q}$ be the two-sided ideal of ${\mathbb K}\langle X,\Theta_1,\Theta_2 \rangle$ which generated by the elements $X\Theta_1-q\Theta_1 X$, $X\Theta_2-p\Theta_2 X$, $\Theta_1\Theta_2+pq^{-2}\Theta_2\Theta_1$, $\Theta_1^2$ and $\Theta_2^2$. The quantum superspace ${\mathbb A}_{p,q}^{1|2}$ with the function algebra
$${\cal O}({\mathbb A}_{p,q}^{1|2}) = {\mathbb K}\langle X,\Theta_1,\Theta_2 \rangle/I_{p,q}$$
is called quantum superspace.
\end{definition}
In accordance with this definition, we have
\begin{equation} \label{6.1}
X\Theta_1 = q\Theta_1 X, \quad X\Theta_2 = p\Theta_2 X, \quad \Theta_1\Theta_2 = -pq^{-2}\Theta_2\Theta_1, \quad \Theta_i^2 = 0
\end{equation}
where $p,q \in {\mathbb K}-\{0\}$.

\begin{example} 
If we consider the generators of the algebra ${\cal O}({\mathbb A}^{1|2}_{p,q})$ as linear maps, then we can find the matrix representations of them. In fact,
it can be seen that there exists a representation $\rho:{\cal O}({\mathbb A}^{1|2}_{p,q})\to \mbox{M}(3,{\mathbb K}')$ such that matrices
\begin{equation*}
\rho(X)=\begin{pmatrix} q & 0 & 0 \\ 0 & p & 0 \\ 0 & 0 & pq \end{pmatrix}, \quad
\rho(\Theta_1)=\begin{pmatrix} 0 & 0 & \varepsilon_1 \\ 0 & 0 & 0 \\ 0 & 0 & 0 \end{pmatrix}, \quad
\rho(\Theta_2)=\begin{pmatrix} 0 & 0 & 0 \\ 0 & 0 & \varepsilon_2 \\ 0 & 0 & 0 \end{pmatrix}
\end{equation*}
representing the coordinate functions satisfy relations $(\ref{6.1})$ for all $\varepsilon_1, \varepsilon_2$.
\end{example}

Let ${\mathbb K}\langle u,\xi_1,\xi_2 \rangle$ be a free algebra generated by $u$, $\xi_1$, $\xi_2$, where $\tau(u)=0$, $\tau(\xi_1)=1=\tau(\xi_2)$. Let ${\cal L}$ be the quotient of the free algebra ${\mathbb K}\langle u,\xi_1,\xi_2 \rangle$ by the two-sided ideal $J_0$ generated by the elements $u\xi_k-\xi_k u$, $\xi_1\xi_2+\xi_2\xi_1$, $\xi_k^2$ for $k=1,2$.

We now define the generators of the algebra ${\cal F}({\mathbb A}_{p,q}^{1|2})$ as
\begin{equation*} \label{6.2}
X := e^u, \quad \Theta_k := e^{ku} \xi_k,
\end{equation*}
for $k=1,2$. Then, by direct calculations we can prove the following lemma.

\begin{lemma} 
The generators $u$, $\xi_1$, $\xi_2$ have the following commutation relations (Lie (anti-)brackets), for $j,k=1,2$
\begin{equation} \label{6.3}
[u,\xi_k] = {\bf i} \,\hbar_k \,\xi_k, \quad [\xi_j,\xi_k]_+ = 0,
\end{equation}
where $q=e^{{\bf i}\,\hbar_1}$, $p=e^{{\bf i}\,\hbar_2}$ with ${\bf i}=\sqrt{-1}$ and $\hbar_1,\hbar_2\in{\mathbb R}$.
\end{lemma}

\noindent
We denote the algebra for which the generators obey the relations (\ref{6.3}) by ${\cal L}_{\hbar_1,\hbar_2}:={\cal L}({\mathbb A}_{p,q}^{1|2})$. The
${\mathbb Z}_2$-graded Hopf algebra structure of ${\cal L}_{\hbar_1,\hbar_2}$ can be read off from Theorem 2.3:

\begin{theorem} 
The Lie superalgebra ${\cal L}_{\hbar_1,\hbar_2}$ is a ${\mathbb Z}_2$-graded Hopf algebra with coproduct, counit and coinverse on the algebra
${\cal L}_{\hbar_1,\hbar_2}$ defined by
\begin{align*}
\Delta(u_i) &= u_i\otimes{\bf 1} + {\bf 1}\otimes u_i, & \epsilon(u_i) &= 0, & S(u_i) &= - u_i.
\end{align*}
for $u_i\in\{u,\xi_1,\xi_2\}$.
\end{theorem}

\begin{example} 
There exists a Lie algebra homomorphism $\mu$ from ${\cal L}_{\hbar_1,\hbar_2}$ into ${\rm M}(3,{\mathbb K}')$.
\end{example}

\noindent{\it Proof} 
We seen that there exists an algebra homomorphism $\rho$ from ${\cal F}({\mathbb A}_{p,q}^{1|2})$ into ${\rm M}(3,{\mathbb K}')$ such that the relations (\ref{6.1}) hold. As a consequence of this fact, there exists a Lie algebra homomorphism $\mu$ from ${\cal L}_{\hbar_1,\hbar_2}$ into
${\rm M}(3,{\mathbb K}')$. The action of $\mu$ on the generators of ${\cal L}_{\hbar_1,\hbar_2}$ is of the form
\begin{align} \label{6.4}
\mu(u) &= \left(\begin{matrix} {\bf i}\hbar_2 & 0 & 0 \\ 0 & {\bf i}\hbar_1 & 0 \\ 0 & 0 & {\bf i}(\hbar_1 + \hbar_2) \end{matrix}\right), \quad
\mu(\xi_1) = \left(\begin{matrix} 0 & 0 & 0 \\ 0 & 0 & 0 \\ e^{-{\bf i}(\hbar_1 + \hbar_2)}\varepsilon_1 & 0 & 0 \end{matrix}\right), \nonumber\\
\mu(\xi_2) &= \left(\begin{matrix} 0 & 0 & 0 \\ 0 & 0 & 0 \\ 0 & e^{-2{\bf i}(\hbar_1 + \hbar_2)}\varepsilon_2 & 0 \end{matrix}\right)
\end{align}
where $\varepsilon_1$ and $\varepsilon_2$ are two Grassmann numbers. To see that the relations (\ref{6.3}) are preserved under the action of $\mu$, we use the fact that
\begin{equation*}
\mu[a,b] = [\mu(a) , \mu(b)],
\end{equation*}
for all $a,b \in {\cal L}_{\hbar_1,\hbar_2}$. $\square$

\section{$\star$-Structures on the algebras ${\cal O}({\mathbb A}_{h}^{1|2})$ and ${\cal O}({\mathbb A}_{h'}^{2|1})$}\label{sec7}
It is possible to define the star operation (or involution) on the Grassmann generators. However, there are two possibilities to do so \cite{FSS}.
If $\alpha$ and $\beta$ are two Grassmann generators and $\lambda$ is a complex number and $\bar{\lambda}$ its complex conjugate, the star operation, denoted by $\star$, is defined by
\begin{equation*}
(\lambda\alpha)^\star = \bar{\lambda}\alpha^\star, \quad (\alpha\beta)^\star = \beta^\star \alpha^\star, \quad (\alpha^\star)^\star = \alpha
\end{equation*}
and the superstar operation, denoted by $\#$, is defined by
\begin{equation*}
(\lambda\alpha)^\# = \bar{\lambda}\alpha^\#, \quad (\alpha\beta)^\# = \alpha^\# \beta^\#, \quad (\alpha^\#)^\# = -\alpha.
\end{equation*}

It is easily shown that there exists a star operation on the algebra ${\cal O}({\mathbb A}_{q}^{1|2})$ if $q$ is a complex number of modulus one:
\begin{proposition} 
(i) If $\bar{q}=q^{-1}$ then the algebra ${\cal O}({\mathbb A}_{q}^{1|2})$ equipped with the involution determined by
\begin{equation} \label{7.1}
X^\star = X, \quad \Theta_i^\star = \Theta_i \qquad (i=1,2)
\end{equation}
becomes a $\star$-algebra.

(ii) If $\bar{p}=p^{-1}$ and $\bar{q}=q^{-1}$ then the algebra ${\cal O}({\mathbb A}_{p,q}^{2|1})$ equipped with the involution determined by
\begin{equation} \label{7.2}
\Phi^\star = \Phi, \quad Y_i^\star = - Y_i \qquad (i=1,2)
\end{equation}
becomes a $\star$-algebra.
\end{proposition}

\subsection{$\star$-Structures on the algebra ${\cal O}({\mathbb A}_{h}^{1|2})$}\label{sec7.1}
As noted in Section 3, the relations in (\ref{3.3}) do not include the parameter $h'$. Thus, we can rearrange the change of basis in the coordinates (see, equation (\ref{3.1})) as
\begin{equation} \label{7.3}
\begin{pmatrix} X \\ \Theta_1 \\ \Theta_2 \end{pmatrix} = \begin{pmatrix} 1 & 0 & 0 \\ 0 & 1 & 0 \\ \frac{h}{q-1} & 0 & 1 \end{pmatrix} \begin{pmatrix} x \\ \theta_1 \\ \theta_2 \end{pmatrix}.
\end{equation}
This case can help us to define a star operation on the algebra ${\cal O}({\mathbb A}_{h}^{1|2})$ by a coordinate transformation using the generators of the algebra ${\cal O}({\mathbb A}_{q}^{1|2})$ and to prove the following lemma.

\begin{lemma} 
For a certain special choice of $h$, there exists an involution on the algebra ${\cal O}({\mathbb A}_{h}^{1|2})$.
\end{lemma}

\noindent{\it Proof} 
Using the equation (\ref{7.3}), we introduce the coordinates $x$, $\theta_1$ and $\theta_2$ with the change of basis in the coordinates of the superspace ${\mathbb A}_{q}^{1|2}$ as follows:
\begin{equation*}
x = X, \quad \theta_1 = \Theta_1, \quad \theta_2 = \Theta_2 - \frac{h}{q-1} \,X.
\end{equation*}
Then, with $|q|=1$ and (\ref{7.1})
\begin{align*}
\theta_2^\star &= \Theta_2^\star - \frac{q\bar{h}}{1-q} \,X^\star = \theta_2 + \frac{h + q\bar{h}}{q-1} \,x
\end{align*}
so that, if we demand that $\bar{h}=-h$, we obtain $\theta_2^\star = \theta_2 - h x$. Note that
\begin{equation*}
(x^\star)^\star = x, \quad (\theta_1^\star)^\star = \theta_1, \quad (\theta_2^\star)^\star = \theta_2,
\end{equation*}
for all $h$. $\square$

\begin{proposition} 
If $\bar{h}=-h$, then the algebra ${\cal O}({\mathbb A}_{h}^{1|2})$ supplied with the involution determined by
\begin{equation} \label{7.4}
x^\star = x, \quad \theta_1^\star = \theta_1, \quad \theta_2^\star = \theta_2 - h x
\end{equation}
becomes a $\star$-algebra.
\end{proposition}

\noindent{\it Proof} 
Since $\bar{h}=-h$, we have
\begin{align*}
(x\theta_1-\theta_1 x)^\star &= \theta_1 x - x\theta_1, \\
(x\theta_2-\theta_2 x - h x^2)^\star &= (\theta_2 - h x) x - x (\theta_2 - h x) + h x^2 = (\theta_2 x - x\theta_2 + h x^2), \\
(\theta_1\theta_2+\theta_2\theta_1)^\star &= (\theta_2 - hx)\theta_1 + \theta_1(\theta_2 - hx) = \theta_2\theta_1 + \theta_1\theta_2,\\
(\theta_2^2 + h \theta_2 x)^\star &= (\theta_2 - hx)(\theta_2 - h x) + x(\theta_2 -h x)(-h) = \theta_2^2 + h \theta_2 x.
\end{align*}
Hence the ideal $(x\theta_1-\theta_1 x, \,x\theta_2-\theta_2 x-hx^2, \,\theta_1\theta_2+\theta_2\theta_1, \,\theta_1^2, \,\theta^2_2+h\theta_2 x)$ is $\star$-invariant and the quotient algebra
\begin{equation*}
{\mathbb K}\langle x,\theta_1,\theta_2\rangle/(x\theta_1-\theta_1 x, \,x\theta_2-\theta_2 x-hx^2, \,\theta_1\theta_2+\theta_2\theta_1, \,\theta_1^2, \,\theta^2_2+h\theta_2 x)
\end{equation*}
becomes a $\star$-algebra. $\square$

\noindent{\bf Note 5.} Of course, we can consider the change of basis in the coordinates of the superspace ${\mathbb A}_q^{1|2}$ in (\ref{3.1}). In this case, since
\begin{align*}
x^\star &= (1 + \tilde{h}\bar{\tilde{h}}' - \overline{\tilde{h}\tilde{h'}})x + (\tilde{h}' - \bar{\tilde{h}}')\theta_2, \\
\theta_1^\star &= \theta_1,\\
\theta_2^\star &= (1 - \bar{\tilde{h}}(\tilde{h}' - \bar{\tilde{h}}'))\theta_2 + (\tilde{h} - \bar{\tilde{h}})x,
\end{align*}
we have again (\ref{7.4}) with the choices $\bar{h}=-h$ and $\bar{h}'=h'$.

\subsection{$\star$-Structure on the algebra ${\cal O}({\mathbb A}_{h'}^{2|1})$}\label{sec7.2}
Since the relations in (\ref{3.6}) do not include the parameter $h$, we can rearrange the change of basis in the coordinates (see, equation (\ref{3.1})) as
\begin{equation} \label{7.5}
\begin{pmatrix} \Phi \\ Y_1 \\ Y_2 \end{pmatrix} = \begin{pmatrix} 1 & 0 & \frac{h'}{pq-1} \\ 0 & 1 & 0 \\ 0 & 0 & 1 \end{pmatrix} \begin{pmatrix} \varphi \\ y_1 \\ y_2 \end{pmatrix}.
\end{equation}
There exists a special case, where the algebra ${\cal O}({\mathbb A}_{h'}^{2|1})$ admits an involution. The proofs of the following lemma and proposition can be done in a similar way to Lemma 7.2 and Proposition 7.3.

\begin{lemma} 
If $\bar{h}'=h'$, there exists an involution on the algebra ${\cal O}({\mathbb A}_{h'}^{2|1})$.
\end{lemma}

\begin{proposition} 
If $\bar{h}'=h'$, then the algebra ${\cal O}({\mathbb A}_{h}^{2|1})$ supplied with the involution determined by
\begin{equation} \label{7.6}
\varphi^\star = \varphi - h'y_2, \quad y_i^\star = -y_i, \quad (i=1,2)
\end{equation}
becomes a $\star$-algebra.
\end{proposition}

\end{document}